\documentclass[10pt]{article} \usepackage{amsfonts} \usepackage{amssymb}\usepackage{amsmath}\usepackage{amsfonts}

\makeatletter
\renewcommand{\@makefnmark}{}
\makeatother
\begin{document}
\baselineskip=10pt
\pagestyle{plain}
{\Large

\newcommand{\om}{O(\frac{1}{\mu})}
\newcommand{\lp}{L_2(0,\pi)}
\newcommand{\lpp}{L_2(\Omega)}

\newcommand{\elx}{e^{i\lambda x}}
\newcommand{\elxx}{e^{-i\lambda x}}

\newcommand{\elp}{e^{i\pi\lambda }}
\newcommand{\elpp}{e^{-i\pi\lambda }}
\newcommand{\No}{\textnumero}

\newcommand{\pp}{\prod_{n=-\infty\atop n\ne0}^\infty}
\newcommand{\sss}{\sum_{n=-\infty}^\infty}
\newcommand{\sN}{\sum_{n=-N}^N}
\newcommand{\sNN}{\sum_{|n|>N}}
\newcommand{\pps}{\sum_{n=-\infty\atop n\ne0}^\infty}
\newcommand{\sNNn}{\sum_{|n|\le N}}
\newcommand{\ppo}{\prod_{n=-\infty}^\infty}
\newcommand{\ppj}{\prod_{n=-\infty\atop n\ne0}^\infty\prod_{j=1}^2}
\newcommand{\ppjj}{\prod_{n=-\infty}^\infty\prod_{j=1}^2}

\newcommand{\ssj}{\sum_{n=-\infty\atop n\ne0}^\infty\sum_{j=1}^2}

\newcommand{\ssjj}{\sum_{n=-\infty}^\infty\sum_{j=1}^2}

\newcommand{\mun}{\mu_{n,j}}
\newcommand{\wnj}{W_{N,j}(\mu)}
\newcommand{\ntm}{|2n-\theta-\mu|}
\newcommand{\pwp}{PW_\pi^-}
\newcommand{\cpm}{\cos\pi\mu}
\newcommand{\spm}{\frac{\sin\pi\mu}{\mu}}
\newcommand{\tet}{(-1)^{\theta+1}}
\newcommand{\dm}{\Delta(\mu)}
\newcommand{\dl}{\Delta(\lambda)}

\newcommand{\dnm}{\Delta_N(\mu)}
\newcommand{\sni}{\sum_{n=N+1}^\infty}
\newcommand{\muk}{\sqrt{\mu^2+q_0}}
\newcommand{\agt}{\alpha, \gamma, \theta,}
\newcommand{\muq}{\sqrt{\mu^2-q_0}}
\newcommand{\smn}{\sum_{n=1}^\infty}
\newcommand{\lop}{{L_2(0,\pi)}}
\newcommand{\tdm}{\tilde\Delta_+(\mu_n)}
\newcommand{\dmn}{\Delta_+(\mu_n)}
\newcommand{\dd}{D_+(\mu_n)}
\newcommand{\ddd}{\tilde D_+(\mu_n)}
\newcommand{\sn}{\sum_{n=1}^\infty}
\newcommand{\pnn}{\prod_{n=1}^\infty}
\newcommand{\aln}{\alpha_n(\mu)}
\newcommand{\emp}{e^{\pi|Im\mu|}}
\newcommand{\ppp}{\prod_{p=p_0}^\infty}
\newcommand{\ppk}{\sum_{p=p_0}^{\infty}\sum_{k=1}^{[\ln p]}}
\newcommand{\mln}{m(\lambda_n)}
\newcommand{\mnk}{\mu_{n_k}}
\newcommand{\sxm}{s(x,\mu)}
\newcommand{\skm}{s(\xi,\mu)}
\newcommand{\cxm}{c(x,\mu)}
\newcommand{\ckm}{c(\xi,\mu)}

\newcommand{\lnk}{\lambda_{n_k}}
\medskip
\medskip
\medskip

\footnote{

\hspace{2mm}Keywords: Dirac operator, non-regular boundary conditions, completeness of root function systems}

\centerline {\bf On the completeness of root function system of the $2\times 2$ Dirac operators}
\centerline {\bf   with non-regular boundary conditions}
\medskip
\medskip
\medskip
\medskip
\medskip

\centerline { Alexander Makin}
\medskip
\medskip
\medskip
{\normalsize
\centerline {Peoples Friendship University of Russia}
 \centerline {	117198, Miklukho-Maklaya str. 6, Moscow, Russia}

\medskip
\medskip
\medskip

\medskip
\medskip
\begin{quote}{\normalsize
The paper is concerned with the completeness property of root functions of  the $2\times 2$ Dirac operator with summable complex-valued
 potential  and non-regular boundary conditions. Sufficient conditions for the completeness of the root function system of the operator under consideration are established.}
   \end{quote}

\centerline {\bf 1. Introduction}

\medskip
\medskip
\medskip
\medskip

In the  present paper, we study the Dirac system
\begin{equation}
B\mathbf{y}'+V\mathbf{y} =\lambda\mathbf{y},\label{1:2}
 \end{equation}
 where $\mathbf{y}={\rm col}(y_1(x),y_2(x))$,
 \[
 B=\begin{pmatrix}
 -i&0\\
 0&i
 \end{pmatrix},\quad V=\begin{pmatrix}
 0&P(x)\\
 Q(x)&0
 \end{pmatrix},
\]
the functions $P, Q\in L_1(0,\pi)$,   with two-point boundary conditions
$$
\begin{array}{c}
U_1(\mathbf{y})= a_{11}y_1(0)+a_{12}y_2(0)+a_{13}y_1(\pi)+ a_{14}y_2(\pi)=0,\\ U_2(\mathbf{y})= a_{21}y_1(0)+a_{22}y_2(0)+a_{23}y_1(\pi)+ a_{24}y_2(\pi)=0,
\end{array}\eqno(2)
$$
where
the coefficients $a_{jk}$ are arbitrary complex numbers,
and rows of  the matrix
\[
A=\begin{pmatrix}
 a_{11}&a_{12}&a_{13} &a_{14}\\
 a_{21}&a_{22}&a_{23} &a_{24}
\end{pmatrix}
\]
are linearly independent.

The operator $\mathbb{L}\mathbf{y}=B\mathbf{y}'+V\mathbf{y}$ is regarded as a linear operator in the space
$\mathbb{H}=L_2(0,\pi)\oplus L_2(0,\pi)$,
with the domain $D(\mathbb{L})=\{\mathbf{y}\in W_1^1[0,\pi]\oplus  W_1^1[0,\pi]:\, \mathbb{L}\mathbf{y}\in \mathbb{H}$, $U_j(\mathbf{y})=0$ $(j=1,2)\}$.

Denote by $A_{jk}$ $(1\le j<k\le4)$ the determinant composed of the jth and kth columns of the matrix $A$.
Boundary conditions (2) are called regular if
$$
A_{14}A_{23}\ne0,
$$
otherwise they are irregular or degenerate.

The general spectral problem for $n\times n$ first order system of ordinary differential equations
(ODE) on a finite interval for the first time has been investigated by G. Birkhoff
and R. Langer \cite{bib13}. More precisely, they introduced the concepts of regular and
strictly regular boundary conditions, investigated the asymptotic behavior of eigenvalues
and eigenfunctions and proved a pointwise convergence result on spectral decompositions
for the corresponding differential operator. The first completeness result for such
systems was established by V.P. Ginzburg \cite{bib14} who treated the case $B = I_n$, $V(\cdot) = 0$.
V.A. Marchenko \cite{bib8} established completeness property for the system of root functions
of the operator $\mathbb{L}$ with regular boundary conditions and continuous matrix potential $V$.
This restriction occurs because the transformation operators used for the proof have
been constructed in \cite{bib8} only for continuous potentials.

Later,  M.M. Malamud and L.L. Oridoroga  \cite{bib7} established completeness property for
$B$-weakly regular boundary value problems for arbitrary $n\times n$ first order systems of
ODE with integrable matrix potential $V\in L^1([0,\pi];\mathbb{C}^{n\times n})$ (originally this result was
announced in \cite{bib16} in 2000).

The first result on completeness for the $2\times 2$  Dirac-type operator  $\mathbb{L}$ with the matrix
$B = diag(b_1, b_2)$ and non-regular boundary conditions was established in \cite{bib7}. Namely,  states that under the smoothness assumption $P,Q\in C^1[0,\pi]$ the system
of root functions of the operator $\mathbb{L}$ is complete whenever both of the following conditions
hold:
$$
\begin{array}{c}
|A_{32}| + |b_1A_{13}P(0) + b_2A_{42}Q(\pi)|\ne0,\\
|A_{14}| + |b_1A_{13}P(\pi) + b_2A_{42}Q(0)|\ne0.
\end{array}\eqno(3)
$$
In \cite{bib1}, similar results were obtained in the case when $B\ne B^*$ and $P,Q$ are analytic.
Both   \cite{bib1} and  \cite{bib7} relied on the method of transformation operators.
In \cite{bib15} and \cite{bib3}, A.A. Lunyov and M.M. Malamud generalized results of \cite{bib7} to
establish potential-dependent completeness and spectral synthesis results for the system
of root functions of the $n\times n$  system with non-weakly-regular boundary conditions
assuming that  $n\times n$  potential matrix $V(\cdot) = 0$ is continuous only at the endpoints $0$ and $\pi$.
In \cite{bib4} the same authors extended completeness results from \cite{bib7} for the $2\times 2$ Dirac-type
operator $\mathbb{L}$  that involve boundary values $V^{(k)}(0)$ and $V^{(k)}(\pi)$, $k\in{0, 1,\ldots, m-1}$, of
the derivatives of the potential $V \in W^m_2([0, \pi];C^{2\times 2})$.
In \cite{bib2} A.P. Kosarev and A.A. Shkalikov extended completeness
results from  \cite{bib3}, \cite{bib7} to the case of $2\times 2$ Dirac-type operators with non-constant matrix
$B = diag(b_1(x), b_2(x))$ and degenerate boundary conditions of a special form $(y_1(0) =
y_2(\pi) = 0)$. Namely, the system of root functions of the operator  $\mathbb{L}$ is complete whenever
the functions $b_1, b_2, P,Q$ are absolutely continuous and satisfy the condition $P(\pi)Q(0)\ne0$.

Notice, if conditions  (2) are not regular the completeness property essentially depends on the potential $V$, in particular, in this case
the root function system of nonperturbed operator

$$
B\mathbf{y}'=\lambda\mathbf{y},\quad U(\mathbf{y})=0\eqno(4)
$$
is not complete in $\mathbb{H}$ \cite{bib8}.

 In a very recent paper \cite{bib77}, A.S. Makin obtained sufficient conditions of
the completeness for the root function system of problem (1), (2) when $A_{14}A_{23}=0$, $|A_{13}|+|A_{42}|>0$ and the potential $V\in L_1(0,\pi)$, and in \cite{bib99} A.A. Lunyov and M.M. Malamud refined their previous results. If the completeness property of a system of root functions is established, the question arises whether it forms a basis.
Most
complete result on the Riesz basis property of boundary value problems for
$2\times2$ Dirac systems with $V(\cdot)\in L^1([0,1]; \mathbb{C}^{2\times2})$$V(\cdot)\in L^1([0,1]; \mathbb{C}^{2\times2})$ and strongly regular
boundary conditions obtained independently and at the same time,
but using different methods by A. M. Savchuk and A. A. Shkalikov \cite{bib33}
on the one hand, and by A. A. Lunyov and M. M. Malamud \cite{bib5}, \cite{bib43} on the other. Block basis
Riesz in the case of an $L_1$-potential matrix and regular boundary
conditions was first proven in \cite{bib33}. In \cite{bib88} the author considered
  spectral problems for the Dirac operator with regular but not strongly regular boundary conditions and complex-valued
   summable  potential and  found  conditions under which the root function system forms a usual Riesz basis rather than a Riesz basis with parentheses.

 In the  present paper, we study the completeness property in the case
 $A_{13}=A_{42}=0$.

\medskip
\medskip

\centerline {\bf 2. Preliminaries}

Denote by
 \[
E(x,\lambda)=\begin{pmatrix}
e_{11}(x,\lambda)&e_{12}(x,\lambda)\\
e_{21}(x,\lambda)&e_{22}(x,\lambda)
\end{pmatrix}\eqno(5)
\]
the matrix of the fundamental solution system  to  equation
(1)
with boundary condition
$
E(0,\lambda)=I
$, where
$I$ is the unit matrix.
It is well known \cite{bib5} that
$$
\begin{array}{c}
e_{11}(x,\lambda)=e^{ix\lambda}(1+o(1))+e^{-ix\lambda}o(1),\quad e_{12}(x,\lambda)=e^{ix\lambda}o(1)+e^{-ix\lambda}o(1),\\

\quad e_{21}(x,\lambda)=e^{ix\lambda}o(1)+e^{-ix\lambda}o(1),\quad e_{22}(x,\lambda)=e^{ix\lambda}o(1)+e^{-ix\lambda}(1+o(1))
\end{array}\eqno(6)
$$
as $\lambda\to\infty$ uniformly in $x\in[0,\pi]$.

The eigenvalues of problem (1), (2) are the roots of the characteristic equation
$$
\Delta(\lambda)=0,
$$
where
$$
\Delta(\lambda)=
\left|\begin{array}{cccc}
U_1(E^{[1]}(\cdot,\lambda))&U_1(E^{[2]}(\cdot,\lambda))\\
U_2(E^{[1]}(\cdot,\lambda))&U_2(E^{[2]}(\cdot,\lambda))\\
\end{array}
\right|,
$$
$E^{[k]}(x,\lambda)$ is the $k$th column of matrix (5).

It was shown in [7] by the method of transformation operators that the characteristic determinant
$\Delta(\lambda)$ of problem (1), (2)
can be reduced to the form
$$
\begin{array}{c}
\Delta(\lambda)=A_{12}+A_{34}+A_{32}e_{11}(\pi,\lambda)+A_{14}e_{22}(\pi,\lambda)+A_{13}e_{12}(\pi,\lambda)+A_{42}e_{21}(\pi,\lambda)=\\\\
=\Delta_0(\lambda)+\int_0^\pi r_1(t)e^{-i\lambda t}dt+\int_0^\pi r_2(t)e^{i\lambda t}dt,
\end{array}\eqno(7)
$$
where the function
$$
\Delta_0(\lambda)=A_{12}+A_{34}-A_{23}e^{i\pi\lambda}+A_{14}e^{-i\pi\lambda}
$$
is the characteristic determinant of problem (4) and the functions $r_j\in L_1(0,\pi)$, $j=1,2$.

Note also that in the recent paper \cite{bib6} representation like (7)  for characteristic determinant was obtained for general first order $n\times n$ -systems of ODE.

 For convenience, we present several commonly used relations established in \cite{bib77}. Let $\lambda$ be a complex number, $Im\lambda\ne0$, $\rho>0$. Suppose $\tau(x)$ is a continuous function on the segment $[0,\pi]$. Then, for any $b\in[0,\pi]$
$$
|\int_0^b x^\rho e^{-2|Im\lambda| x}\tau(x)dx|\le\frac{c}{|Im\lambda|^{\rho+1}},\eqno(8)
$$
where $c$  not depending on $b$.
If a function $\tau\in L_1(0,\pi)$ then the following relation is valid
$$
|\int_0^\pi x^\rho e^{-2|Im\lambda| x}\tau(x)dx|=\frac{o(1)}{|Im\lambda|^{\rho}}
\eqno(9)
$$
as $|Im\mu|\to\infty$.
In addition, simple computations show that if $\rho>0, \lambda>0$,   $\rho\le \pi\lambda$, then
$$
\max_{0\le x\le \pi}x^\rho e^{-\lambda x}=\frac{\rho^\rho}{\lambda^\rho}e^{-\rho}. \eqno(10)
$$

Denote
$$
g_0(t,\lambda)=1,\eqno(11)
$$

$$
g_1(t,\lambda)=\int_0^t e^{-2i\lambda t_1}P(t_1)dt_1\int_0^{t_1}e^{2i\lambda t_2}Q(t_2)dt_2,\eqno(12)
$$
$$
\begin{array}{c}
g_n(t,\lambda)=\int_0^t e^{-2i\lambda t_1}P(t_1)dt_1\int_0^{t_1}e^{2i\lambda t_2}Q(t_2)dt_2\ldots\\\\\ldots
\int_0^{t_{2n-2}}e^{-2i\lambda t_{2n-1}}P(t_{2n-1})dt_{2n-1}\int_0^{t_{2n-1}}e^{2i\lambda t_{2n}}Q(t_{2n})dt_{2n}
\end{array}\eqno(13)
$$
and, analogously,
denote
$$
h_0(t,\lambda)=1,\eqno(14)
$$

$$
h_1(t,\lambda)=\int_0^t e^{2i\lambda t_1}Q(t_1)dt_1\int_0^{t_1}e^{-2i\lambda t_2}P(t_2)dt_2,\eqno(15)
$$
$$
\begin{array}{c}
h_n(t,\lambda)=\int_0^t e^{2i\lambda t_1}Q(t_1)dt_1\int_0^{t_1}e^{-2i\lambda t_2}P(t_2)dt_2\ldots\\\\\ldots
\int_0^{t_{2n-2}} e^{2i\lambda t_{2n-1}}Q(t_{2n-1})dt_{2n-1}\int_0^{t_{2n-1}}e^{-2i\lambda t_{2n}}P(t_2)dt_{2n}.
\end{array}\eqno(16)
$$

{\bf Lemma 1.}\cite{bib77} {\it The following representations are valid

$$
e_{11}(t,\lambda)=e^{i\lambda t}\sum_{n=0}^\infty g_n(t,\lambda),\eqno(17)
$$
$$
e_{22}(t,\lambda)=e^{-i\lambda t}\sum_{n=0}^\infty h_n(t,\lambda),\eqno(18)
$$
where the series in right-hand sides of (17-18) for any $\lambda$ converge uniformly and absolutely on the segment $[0,\pi]$. }

\medskip
\medskip
\medskip
\centerline {\bf 3. Main results}
\medskip
\medskip
\medskip

Let $0<\varepsilon<\pi/10$.
Denote by $\Omega^+_\varepsilon$ the domain  $\varepsilon\le \arg\lambda\le\pi-\varepsilon$, and by $\Omega^-_\varepsilon$ the domain
 $-\pi+\varepsilon\le \arg\lambda\le-\varepsilon$. Further, $\|f\|$ stands for $\|f\|_{L_1(0,\pi)}$.

{\bf Lemma 2.} {\it Suppose

$$
\lim_{h\to0}\frac{\int_{\pi-h}^{\pi}P(x)dx}{h^{\rho_4}}=\nu_4\ne0,\quad\lim_{h\to0}\frac{\int_0^hQ(x)dx}{h^{\rho_6}}=\nu_6\ne0,
\eqno(19)
$$
where $\rho_4>0$, $\rho_6>0$.

Then in the domain $\Omega_\varepsilon^+$
$$
|e_{11}(\pi,\lambda)|\ge\frac{c_1e^{\pi|{\rm Im}\lambda|}}{|{\rm Im}\lambda|^{\rho_4+\rho_6}},\eqno(20)
$$

where $c_1>0$.

}
Proof. First of all, we establish a number of inequalities that will be used later. Let $0\le t\le\pi$.
Integrating by parts, we obtain
$$
\begin{array}{c}
 \int_0^t e^{2i\lambda y}Q(y)dy=
 e^{2it\lambda}\int_0^t Q(y)dy-2i\lambda \int_0^t e^{2i\lambda y}dy(\int_0^y
 Q(t_1)dt_1)=\\\\=
 t^{\rho_6} e^{2it\lambda}(\nu_6 + \tau(t))-2i\lambda \int_0^t y^{\rho_6} e^{2i\lambda y}(\nu_6 + \tau(y))dy.
 \end{array}
$$
This together with (8)  and (10) implies

$$
|\int_0^{t} e^{2i\lambda x}Q(x)dx|\le\frac{c_2}{|Im\lambda|^{\rho_6}},\eqno(21)
$$
where $c_2$ does not depend on $t$.

Denote
$\hat P(x)= P(\pi-x)$. Then integrating by parts and replacing $\pi-t=v$, we obtain

$$
\begin{array}{c}
 \int_t^\pi e^{-2i\lambda y}P(y)dy= e^{-2i\pi\lambda}\int_0^{\pi-t}e^{2is\lambda}\hat P(s)ds=\\\\=
 e^{-2it\lambda}\int_0^{\pi-t}\hat P(x)dx-2i\lambda e^{-2i\pi\lambda} \int_0^{\pi-t} e^{2i\lambda s}ds(\int_0^s\hat P(x)ds)=\\\\=
(\pi-t)^{\rho_4} e^{-2it\lambda}(\nu_4 + \tau(\pi-t))-2i\lambda e^{-2i\pi\lambda}\int_0^{\pi-t} s^{\rho_4} e^{2i\lambda s}(\nu_4 + \tau(s))ds=\\\\=
 e^{-2i\pi\lambda}(v^{\rho_4} e^{2iv\lambda}(\nu_4 + \tau(v))-2i\lambda\int_0^{\pi-t} s^{\rho_4} e^{2i\lambda s}(\nu_4 + \tau(s))ds.
\end{array}
$$
This together with (8)  and (10) implies

$$
|\int_{t}^\pi e^{-2i\lambda t_1}P(t_1)dt_1|\le\frac{c_3e^{2\pi|{\rm Im}\lambda|}}{|{\rm Im}\lambda|^{\rho_4}}\eqno(22)
$$

and

$$
|\int_0^{\pi-t}e^{2i\lambda s}\hat P(s)ds|\le\frac{c_3}{|{\rm Im}\lambda|^{\rho_4}},\eqno(23)
$$
where in both cases $c_3$ does not depend on $t$.

 Let us estimate the function
$$
g_1(\pi,\lambda)=\int_0^\pi e^{-2i\lambda t}P(t)dt\int_0^{t}e^{2i\lambda x}Q(x)dx.
$$
Changing the order of integration we obtain
$$
\begin{array}{c}
g_1(\pi,\lambda)=\int_0^\pi e^{2i\lambda x}Q(x)dx\int_x^\pi e^{-2i\lambda t}P(t)dt=\\\\=
\int_0^\pi e^{2i\lambda x}Q(x)dx(\int_0^\pi e^{-2i\lambda t}P(t)dt-\int_0^x e^{-2i\lambda t}P(t)dt)=\\\\=
\int_0^\pi e^{2i\lambda x}Q(x)dx\int_0^\pi e^{-2i\lambda t}P(t)dt-\int_0^\pi e^{2i\lambda x}Q(x)dx\int_0^x e^{-2i\lambda t}P(t)dt.
\end{array}
$$
It follows from \cite[ Lemma 3.5]{bib77} that
$$
|\int_0^\pi e^{-2i\lambda t}P(t)dt|\ge\frac{c_4e^{2\pi|{\rm Im}\lambda|}}{|{\rm Im}\lambda|^{\rho_4}}\eqno(24)
$$
$(c_4>0)$.
It follows from \cite[ Lemma 3.7]{bib77} that
$$
|\int_0^\pi e^{2i\lambda x}Q(x)dx|\ge\frac{c_5}{|{\rm Im}\lambda|^{\rho_6}}\eqno(25)
$$
$(c_5>0)$.
Invoking the Holder inequality we have
$$
\begin{array}{c}
|\int_0^\pi e^{2i\lambda x}Q(x)dx\int_0^x e^{-2i\lambda t}P(t)dt|\le\int_0^\pi |Q(x)|dx\int_0^x |e^{2i\lambda (x-t)}|P(t)|dt\le
\|Q\|\|P\|<c_6.
\end{array}
$$
This together with (24) and (25) imply (26)

$$
|g_1(\pi,\lambda)|\ge\frac{c_7e^{2\pi|{\rm Im}\lambda|}}{|{\rm Im}\lambda|^{\rho_4+\rho_6}}\eqno(26)
$$
$(c_7>0)$.

 Let us estimate the function $g_2(\pi,\lambda)$. Using relations (12), (13) and changing the order of integration we obtain

$$
\begin{array}{c}
g_2(\pi,\lambda)=\int_0^\pi e^{-2i\lambda t_1}P(t_1)dt_1\int_0^{t_1}e^{2i\lambda t_2}Q(t_2)g_1(t_2,\lambda)dt_2=\\\\=
\int_0^\pi e^{2i\lambda t_2}Q(t_2)g_1(t_2,\lambda)dt_2\int_{t_2}^\pi e^{-2i\lambda t_1}P(t_1)dt_1
\end{array}\eqno(27)
$$
and
$$
\begin{array}{c}
g_1(t_2,\lambda)=\int_0^{t_2} e^{2i\lambda x}Q(x)dx\int_x^{t_2} e^{-2i\lambda t}P(t)dt=\\\\=
\int_0^{t_2} e^{2i\lambda x}Q(x)dx(\int_0^{t_2} e^{-2i\lambda t}P(t)dt-\int_0^x e^{-2i\lambda t}P(t)dt)=\\\\=
\int_0^{t_2} e^{2i\lambda x}Q(x)dx\int_0^{t_2} e^{-2i\lambda t}P(t)dt-\int_0^{t_2} e^{2i\lambda x}Q(x)dx\int_0^x e^{-2i\lambda t}P(t)dt,
\end{array}\eqno(28)
$$
hence,

$$
\begin{array}{c}
g_2(\pi,\lambda)=
\int_0^\pi e^{2i\lambda t_2}Q(t_2)dt_2\int_0^{t_2} e^{2i\lambda x}Q(x)dx\int_0^{t_2} e^{-2i\lambda t}P(t)dt\int_{t_2}^\pi e^{-2i\lambda t_1}P(t_1)dt_1-\\\\-\int_0^\pi e^{2i\lambda t_2}Q(t_2)dt_2\int_0^{t_2} e^{2i\lambda x}Q(x)dx\int_0^x e^{-2i\lambda t}P(t)dt\int_{t_2}^\pi e^{-2i\lambda t_1}P(t_1)dt_1=I_1-I_2.
\end{array}
$$

It follows from \cite[ Lemma 3.4]{bib77} that
$$
|\int_0^{t_2} e^{-2i\lambda t}P(t)dt|=o(1)e^{2t_2|{\rm Im}\lambda|}.\eqno(29)
$$

The Holder inequality and estimates (21), (22), (29) imply
$$
\begin{array}{c}
|I_1|\le
||Q||\max_{0\le t_2\le\pi}|e^{2i\lambda t_2}\int_0^{t_2} e^{2i\lambda x}Q(x)dx\int_0^{t_2} e^{-2i\lambda t}P(t)dt\int_{t_2}^\pi e^{-2i\lambda t_1}P(t_1)dt_1|=\\\\=\frac{e^{2\pi|{\rm Im}\lambda|}o(1)}{|{\rm Im}\lambda|^{\rho_4+\rho_6}}.
\end{array}\eqno(30)
$$

Consider the addend $I_2$.
Obviously,
$$
I_2=\int_0^\pi e^{2i\lambda t_2}Q(t_2)dt_2\psi(t_2,\lambda)\int_{t_2}^\pi e^{-2i\lambda t_1}P(t_1)dt_1,\eqno(31)
$$
where
$$
\psi(t_2,\lambda)=\int_0^{t_2} e^{2i\lambda x}Q(x)dx\int_0^x e^{-2i\lambda t}P(t)dt.
$$
Condition (19) implies
$$
\int_0^{x}Q(s)ds=\gamma_6x^{\rho_6}+\gamma_6x^{\rho_6}\tau(x),\eqno(32)
$$
where the function $\tau(x)$ is continuous on $[0,\pi]$ and $\tau(0)=0$.
Integrating by parts and using (32) we obtain

$$
\begin{array}{c}
\psi(t_2,\lambda)=\int_0^{t_2}[ e^{2i\lambda x}\int_0^x e^{-2i\lambda t}P(t)dt]d\int_0^xQ(s)ds=\\\\=
e^{2i\lambda t_2}\int_0^{t_2}e^{-2i\lambda t}P(t)dt\int_0^{t_2}Q(s)ds-\\\\-
2i\lambda\int_0^{t_2} e^{2i\lambda x}dx\int_0^x e^{-2i\lambda t}P(t)dt\int_0^xQ(s)ds-\int_0^{t_2}P(x)dx\int_0^xQ(s)ds=\\\\=
[\gamma_6t_2^{\rho_6}+\gamma_6t_2^{\rho_6}\tau(t_2)]e^{2i\lambda t_2}\int_0^{t_2}e^{-2i\lambda t}P(t)dt-\\\\-
2i\lambda\int_0^{t_2} [\gamma_6x^{\rho_6}+\gamma_6x^{\rho_6}\tau(x)]e^{2i\lambda x}dx\int_0^x e^{-2i\lambda t}P(t)dt-
\\\\-\int_0^{t_2}[\gamma_6x^{\rho_6}+\gamma_6x^{\rho_6}\tau(x)]P(x)dx.
\end{array}\eqno(33)
$$
Substituting (33) into (31) we have

$$
\begin{array}{c}
\int_0^\pi e^{2i\lambda t_2}Q(t_2)dt_2\{[\gamma_6t_2^{\rho_6}+\gamma_6t_2^{\rho_6}\tau(t_2)]e^{2i\lambda t_2}\int_0^{t_2}e^{-2i\lambda t}P(t)dt-\\\\-
2i\lambda\int_0^{t_2} [\gamma_6x^{\rho_6}+\gamma_6x^{\rho_6}\tau(x)]e^{2i\lambda x}dx\int_0^x e^{-2i\lambda t}P(t)dt-\\
\\-\int_0^{t_2}[\gamma_6x^{\rho_6}+\gamma_6x^{\rho_6}\tau(x)]P(x)dx\}\int_{t_2}^\pi e^{-2i\lambda t_1}P(t_1)dt_1=\\\\\\
=\gamma_6\{\int_0^\pi [t_2^{\rho_6}+t_2^{\rho_6}\tau(t_2)]e^{2i\lambda t_2}Q(t_2)dt_2e^{2i\lambda t_2}\int_0^{t_2}e^{-2i\lambda t}P(t)dt\int_{t_2}^\pi e^{-2i\lambda t_1}P(t_1)dt_1-\\\\-2i\lambda
\int_0^\pi e^{2i\lambda t_2}Q(t_2)dt_2\int_0^{t_2} [x^{\rho_6}+x^{\rho_6}\tau(x)]e^{2i\lambda x}dx\int_0^x e^{-2i\lambda t}P(t)dt\int_{t_2}^\pi e^{-2i\lambda t_1}P(t_1)dt_1-\\\\-
\int_0^\pi e^{2i\lambda t_2}Q(t_2)dt_2\int_0^{t_2}[x^{\rho_6}+x^{\rho_6}\tau(x)]P(x)dx\int_{t_2}^\pi e^{-2i\lambda t_1}P(t_1)dt_1\}=\\\\=\gamma_6\{I_{21}-I_{22}-I_{23}\}.
\end{array}\eqno(34)
$$
Using the Holder inequality and inequalities (10), (22), (29) we obtain

$$
\begin{array}{c}
|I_{21}|\le||Q||\max_{0\le t_2\le\pi}|t_2^{\rho_6}e^{2i\lambda t_2}||\max_{0\le t_2\le\pi}|e^{2i\lambda t_2}\int_0^{t_2}e^{-2i\lambda t}P(t)dt|\times\\\\\times\max_{0\le t_2\le\pi}|\int_{t_2}^\pi e^{-2i\lambda t_1}P(t_1)dt_1|=
\frac{e^{2\pi|{\rm Im}\lambda|}o(1)}{|{\rm Im}\lambda|^{\rho_4+\rho_6}}.
\end{array}\eqno(35)
$$
It follows from the Holder inequality and inequalities (8), $x\le t_2$,  \cite[ Lemma 3.4]{bib77}, and (22) that
$$
\begin{array}{c}
|I_{22}|\le||Q|||\lambda|\times\\\\\times\max_{0\le t_2\le\pi}| e^{2i\lambda t_2}\int_0^{t_2} [x^{\rho_6}+x^{\rho_6}\tau(x)]e^{2i\lambda x}dx\int_0^x e^{-2i\lambda t}P(t)dt\int_{t_2}^\pi e^{-2i\lambda t_1}P(t_1)dt_1|\le\\\\
\le|Q|||\lambda|\max_{0\le t_2\le\pi}|\int_0^{t_2} |x^{\rho_6}+x^{\rho_6}\tau(x)|e^{-2|Im\lambda| x}dx|\times\\\\\times\max_{0\le t_2\le\pi}|e^{2i\lambda t_2}\int_0^x e^{-2i\lambda t}P(t)dt||\max_{0\le t_2\le\pi}|\int_{t_2}^\pi e^{-2i\lambda t_1}P(t_1)dt_1|=\frac{e^{2\pi|{\rm Im}\lambda|}o(1)}{|{\rm Im}\lambda|^{\rho_4+\rho_6}}.
\end{array}\eqno(36)
$$
It follows from   the Holder inequality, (9), (22) that
$$
\begin{array}{c}
|I_{23}|\le c_8|\int_0^\pi t_2^{\rho_6}e^{-2|Im\lambda| t_2}|Q(t_2)|dt_2\int_0^{t_2}|P(x)|dx|\int_{t_2}^\pi e^{-2i\lambda t_1}P(t_1)dt_1|\le\\\\\le
c_9\|P\|\max_{0\le t_2\le\pi}|\int_{t_2}^\pi e^{-2i\lambda t_1}P(t_1)dt_1|\int_0^\pi t_2^{\rho_6}e^{-2|{\rm Im}\lambda| t_2}|Q(t_2)|dt_2=\\\\
=\frac{e^{2\pi|{\rm Im}\lambda|}o(1)}{|{\rm Im}\lambda|^{\rho_4+\rho_6}}.
\end{array}\eqno(37)
$$
This together with (34) and (36) implies
$$
|I_{2}|=\frac{e^{2\pi|Im\lambda|}o(1)}{|{\rm Im}\lambda|^{\rho_4+\rho_6}}.
$$
Combining the last inequality and (30), we have
$$
g_2(\pi,\lambda)=\frac{e^{2\pi|{\rm Im}\lambda|}o(1)}{|{\rm Im}\lambda|^{\rho_4+\rho_6}}.\eqno(38)
$$

Suppose $n>2$.
Denote
 $$
\begin{array}{c}
F_n(t_2,\lambda)=\int_0^{t_2} e^{-2i\lambda t_3}P(t_3)dt_3\int_0^{t_3}e^{2i\lambda t_4}Q(t_4)dt_4\ldots\\\\\ldots
\int_0^{t_{2n-2}}e^{-2i\lambda t_{2n-1}}P(t_{2n-1})dt_{2n-1}\int_0^{t_{2n-1}}e^{2i\lambda t_{2n}}Q(t_{2n})dt_{2n}.
\end{array}
$$

It is easy to see that
$$
\begin{array}{c}
g_n(\pi,\lambda)=\int_0^\pi e^{-2i\lambda t_1}P(t_1)dt_1\int_0^{t_1}e^{2i\lambda t_2}Q(t_2)F_n(t_2,\lambda)dt_2=\\\\=

\int_0^\pi e^{-2i\lambda t_1}P(t_1)dt_1(\int_0^\pi e^{2i\lambda t_2}Q(t_2)F_n(t_2,\lambda)dt_2-
\int_{t_1}^\pi e^{2i\lambda t_2}Q(t_2)F_n(t_2,\lambda)dt_2)=\\\\=
q_n(\pi,\lambda)\int_0^\pi e^{-2i\lambda t_1}P(t_1)dt_1-
\int_0^\pi e^{-2i\lambda t_1}P(t_1)dt_1
\int_{t_1}^\pi e^{2i\lambda t_2}Q(t_2)F_n(t_2,\lambda)dt_2,
\end{array}\eqno(39)
$$
where

$$
\begin{array}{c}
q_n(t,\lambda)=\int_0^t e^{2i\lambda t_2}Q(t_2)F_n(t_2,\lambda)dt_2=\\\\=
\int_0^t e^{2i\lambda t_2}Q(t_2)dt_2\int_0^{t_2}e^{-2i\lambda t_3}P(t_3)dt_3\ldots\int_0^{t_{2n-1}} e^{2i\lambda t_{2n}}Q(t_{2n})dt_{2n}.
\end{array}\eqno(40)
$$
Consider the first addend in right-hand side of (39). It follows from \cite[ Lemma 3.5]{bib77} that
$$
|\int_0^\pi e^{-2i\lambda t_1}P(t_1)dt_1|\le\frac{c_{10}e^{2\pi|{\rm Im}\lambda|}}{|{\rm Im}\lambda|^{\rho_4}}.
$$

It follows from \cite[ Lemma 3.7]{bib77} that
$$
\sum_{n=3}^\infty |q_n(\pi,\lambda)|=\frac{o(1)}{|{\rm Im}\lambda|^{\rho_6}},
$$
hence,
$$
\sum_{n=3}^\infty |q_n(\pi,\lambda)||\int_0^\pi e^{-2i\lambda t_1}P(t_1)dt_1|=\frac{e^{2\pi|{\rm Im}\lambda|}o(1)}{|{\rm Im}\lambda|^{\rho_4+\rho_6}}.\eqno(41)
$$

Let us estimate the second addend in (39). Changing the order of integration and replacing $t_1=\pi-s$ we obtain

$$
\begin{array}{c}
\int_0^\pi e^{-2i\lambda t_1}P(t_1)dt_1
\int_{t_1}^\pi e^{2i\lambda t_2}Q(t_2)F_n(t_2,\lambda)dt_2=\\\\=
\int_0^\pi e^{2i\lambda t_2}Q(t_2)F_n(t_2,\lambda)dt_2\int_0^{t_2}e^{-2i\lambda t_1}P(t_1)dt_1=\\\\=
e^{-2i\pi\lambda}\int_0^\pi e^{2i\lambda t_2}Q(t_2)F_n(t_2,\lambda)dt_2\int_0^{\pi-t_2}e^{2i\lambda s}\hat P(s)ds,
\end{array}
$$
where $\hat P(s)=P(\pi-s)$.
The Holder inequality and (23) imply
$$
|I|\le\frac{c_{11}e^{2\pi|{\rm Im}\lambda|}}{|{\rm Im}\lambda|^{\rho_4}}\|Q\|\max_{0\le t_2\le\pi}|e^{2i\lambda t_2}F_n(t_2,\lambda)|.\eqno(42)
$$

Consider the function $F_n(t_2,\lambda)$. Denote

$$
\begin{array}{c}
\phi(t_4,\lambda)=e^{2i\lambda t_4}Q(t_4)\int_0^{t_{4}}e^{-2i\lambda t_{5}}P(t_{5})dt_{5}\ldots\\\\\ldots
\int_0^{t_{2n-2}}e^{-2i\lambda t_{2n-1}}P(t_{2n-1})dt_{2n-1}\int_0^{t_{2n-1}}e^{2i\lambda t_{2n}}Q(t_{2n})dt_{2n}.
\end{array}\eqno(43)
$$
Changing the order of integration we obtain

$$
\begin{array}{c}
F_n(t_2,\lambda)=\int_0^{t_2} e^{-2i\lambda t_3}P(t_3)dt_3\int_0^{t_3}\phi(t_4,\lambda)dt_4=\\\\=
\int_0^{t_2}\phi(t_4,\lambda)dt_4\int_{t_4}^{t_2}e^{-2i\lambda t_3}P(t_3)dt_3,
\end{array}\eqno(44)
$$
hence,
$$
|F_n(t_2,\lambda)|\le\int_0^{t_2}|\phi(t_4,\lambda)|dt_4|\int_{t_4}^{t_2}e^{-2i\lambda t_3}P(t_3)dt_3|.\eqno(45)
$$
It follows from \cite[ Lemma 3.4]{bib77} that

$$
|\int_{t_4}^{t_2}e^{-2i\lambda t_3}P(t_3)dt_3|=o(1)e^{2t_2|{\rm Im}\lambda|}.\eqno(46)
$$
Relations (45) and (46) imply
$$
|e^{2i\lambda t_2}F_n(t_2,\lambda)|=o(1)\int_0^{t_2}|\phi(t_4,\lambda)|dt_4.\eqno(47)
$$
Let us estimate the integral in   right-hand side of (47).

It follows from \cite[ Lemma 3.1]{bib77}, (21) and inequality $t_j\ge t_{j+1}$ that
$$
\begin{array}{c}
\int_0^{t_2}|\phi(t_4,\lambda)|dt_4\le\|Q\|\max_{0\le t_4\le\pi}|e^{2i\lambda t_4}\int_0^{t_{4}}e^{-2i\lambda t_{5}}P(t_{5})dt_{5}\ldots\\\\\ldots
\int_0^{t_{2n-2}}e^{-2i\lambda t_{2n-1}}P(t_{2n-1})dt_{2n-1}\int_0^{t_{2n-1}}e^{2i\lambda t_{2n}}Q(t_{2n})dt_{2n}|\le\\\\\le
\|Q\|\max_{0\le t_4\le\pi}|\int_0^{t_{4}}|P(t_{5})|dt_{5}\ldots\\\\\ldots
\int_0^{t_{2n-2}}|e^{2i\lambda(t_4-t_5+t_6\ldots-t_{2n-1})}||P(t_{2n-1})|dt_{2n-1}|\int_0^{t_{2n-1}}e^{2i\lambda t_{2n}}Q(t_{2n})dt_{2n}|\le
\frac{c_{12}^n}{(2n-4)!|{\rm Im}\lambda|^{\rho_6}}.
\end{array}\eqno(48)
$$
Combining (42), (47), (48) we have
$$
|I|=\frac{e^{2\pi|Im\lambda|}o(1)}{|{\rm Im}\lambda|^{\rho_4+\rho_6}}.\eqno(49)
$$
It follows from (39), (41), (49) that
$$
\sum_{n=3}^\infty |g_n(\pi,\lambda)|=\frac{e^{2\pi|{\rm Im}\lambda|}o(1)}{|{\rm Im}\lambda|^{\rho_4+\rho_6}}.\eqno(50)
$$

Combining (17), (38), (26), (50)  we obtain (20).

{\bf Lemma 3.} {\it Suppose

 $$
 A_{14}\ne0,\quad\lim_{h\to0}\frac{\int_0^hP(x)dx}{h^{\rho_5}}=\nu_5\ne0,\quad\lim_{h\to0}\frac{\int_{\pi-h}^{\pi}Q(x)dx}{h^{\rho_7}}=\nu_7\ne0,
 \eqno(51)
 $$
where $\rho_5>0$, $\rho_7>0$;

Then in the domain $\Omega_\varepsilon^-$
$$
|e_{22}(\pi,\lambda)|\ge\frac{ce^{\pi|{\rm Im}\lambda|}}{|{\rm Im}\lambda|^{\rho_5+\rho_7}},\eqno(52)
$$

where $c>0$.

}
Proof. Reasoning as above, it is easy to prove the Lemma 3. Indeed,
from (11-18) it follows that if in formula (17) for the function $e_{11}(\cdot,\cdot)$ we replace $\lambda$ by  $-\lambda$, swap the functions $P$ and $Q$, then we get the function $e_{22}(\cdot,\cdot)$.

Our main result is the following.

{\bf Theorem 1.} {\it
Suppose $A_{14}A_{32}=A_{13}=A_{24}=0$ and
   one of  conditions (53), (54) is satisfied

 $$
 A_{14}\ne0,\quad\lim_{h\to0}\frac{\int_0^hP(x)dx}{h^{\rho_5}}=\nu_5\ne0,\quad\lim_{h\to0}\frac{\int_{\pi-h}^{\pi}Q(x)dx}{h^{\rho_7}}=\nu_7\ne0,
 \eqno(53)
 $$
where $\rho_5>0$, $\rho_7>0$;

$$
A_{32}\ne0,\quad\lim_{h\to0}\frac{\int_{\pi-h}^{\pi}P(x)dx}{h^{\rho_4}}=\nu_4\ne0,\quad\lim_{h\to0}\frac{\int_0^hQ(x)dx}{h^{\rho_6}}=\nu_6\ne0,
\eqno(54)
$$
where $\rho_4>0$, $\rho_6>0$.
Then, the root function system of problem (1), (2) is complete and minimal in $\mathbb{H}$.

}

Proof.
Let $|\lambda|$ be sufficiently large. If condition (53) holds, then it follows from (6) that

$$
|e_{22}(\pi,\lambda)|\ge c_1e^{\pi|{\rm Im} \lambda|}.
$$
if $\lambda\in\Omega_\varepsilon^+$.
It follows from Lemma 3 that
$$
|e_{22}(\pi,\lambda)|\ge\frac{c_2e^{\pi|{\rm Im}\lambda|}}{|{\rm Im}\lambda|^{\rho_5+\rho_7}}\eqno(55)
$$
if $\lambda\in\Omega_\varepsilon^-$
This together with (7) implies that in the domain $\Omega_\varepsilon=\Omega_\varepsilon^-\bigcup\Omega_\varepsilon^+$
$$
|\Delta(\lambda)|\ge \frac{c_3e^{\pi|{\rm Im}\lambda|}}{|{\rm Im}\lambda|^{\rho_5+\rho_7}}.\eqno(56)
$$
Reasoning as above, one can prove the inequality
$$
|\Delta(\lambda)|\ge \frac{c_4e^{\pi|{\rm Im}\lambda|}}{|{\rm Im}\lambda|^{\rho_4+\rho_6}}\eqno(57)
$$
if condition (54) holds and $\lambda\in\Omega_\varepsilon$,
hence, by \cite[ Th. 2.3]{bib3} in both cases the root function system of problem (1), (2) is complete and minimal in $\mathbb{H}$.

\medskip
\medskip
\medskip
\medskip
\medskip

email: alexmakin@yandex.ru

\end{document}